\documentclass{article}
\usepackage{amssymb}
\usepackage{amsmath}
\usepackage{amsfonts,amsthm}

\newtheorem{theorem}{Theorem}
\newtheorem*{acknowledgement}{Acknowledgement}

\newtheorem{definition}[theorem]{Definition}
\newtheorem{example}[theorem]{Example}
\newtheorem{lemma}[theorem]{Lemma}

\newtheorem{proposition}[theorem]{Proposition}
\newtheorem*{remark}{Remark}
\newtheorem{question}[theorem]{Question}
\newcommand{\Sym}{\mathrm{Sym}}
\newcommand{\Alt}{\mathrm{Alt}}
\newcommand{\Aut}{\mathrm{Aut}}
\newcommand{\SL}{\mathrm{SL}}
\newcommand{\EL}{\mathrm{EL}}
\newcommand{\F}{\mathbb{F}}
\newcommand{\G}{\mathfrak{G}}
\begin{document}
\renewcommand{\thefootnote}{}
\title{Cartesian products as profinite completions}
\author{Martin Kassabov \\ Nikolay Nikolov}
 \date{}
\maketitle

\begin{abstract} 
We prove that if a Cartesian product of alternating groups is topologically
finitely generated, then it is the profinite completion of a finitely generated 
residually finite group. 
The same holds for Cartesian producs of other simple groups under some natural restrictions.
\end{abstract}
\section{Introduction}
Given a profinite group $\mathfrak{G}$ it is generally a difficult question to 
determine if there exists a finitely generated residually finite group $G$ 
such that $\mathfrak{G}$ is isomorphic to the profinite completion $\widehat G$ of $G$. 
If this is the case, we then say that $\mathfrak{G}$ is a profinite completion.

Of course for this to happen it is necessary that $\mathfrak{G}$ is topologically
finitely generated but not much is known beyond that. Dan Segal~\cite{dan} has
proved that any collection
of nonabelian finite simple groups can appear as the 
upper composition factors of such a $\mathfrak{G} \simeq \widehat G$.
In his examples the group $G$ is a branch group, i.e., a subgroup of the
automorphism group of a rooted tree with some nice
'branching' properties.

Not all finitely generated profinite groups are profinite completions. 
The argument in~\cite{profinite} 6.2 shows that:

\begin{example} \label{ex} For each $d \in \mathbb N$ and
every infinite sequence of primes $p_1,p_2, \ldots $
the $2$-generated profinite group $\prod_{i=1}^\infty \SL_d(\F_{p_i})$ 
is not a profinite completion.
\end{example}

On the other hand Laci Pyber~\cite{laci} has found examples of 
profinite completions of the form
$$
 \mathfrak G= \widehat {\mathbb{Z}} \times \prod_{n=5}^\infty \Alt(n)^{f(n)}
$$
where the sequence $\{f(n)\}$ satisfies some mild growth conditions,
in particular $f(n)$ can be all $1$ or all $n!$.

It is the nature of these examples that they always have a direct factor
$\hat{\mathbb Z}$ and one might wonder if this is inevitable. 
This was disproved by M.~Kassabov using ideas from~\cite{martin}
and~\cite{martin2} who found a dense finitely generated subgroup of 
$\prod_{n=5}^{\infty}\Alt (n)$ which has property $\tau$ and
with profinite completion
$$
 (\widehat{G})^6 \times \prod_{n=5}^{\infty}\Alt (n) , \quad G=\SL_3 (\F_p[t,t^{-1}]).
$$

In this note we show that for
a Cartesian product of alternating groups the obvious necessary condition 
for being a profinite completion
is also sufficient, i.e.:
\begin{theorem}
\label{main}
If a Cartesian product
$\mathfrak{G} =\prod_{n=5}^\infty \Alt(n)^{f(n)}$ of alternating groups is topologically
finitely generated then it is a profinite completion.
\end{theorem}

We obtain this as a corollary of the following more general result. First, a definition:

\begin{definition}
For a finite simple group $S$ define $l(S)$ to be the 
largest integer $l$ such that $S$ contains a copy of $\Alt(l)$.
 
\end{definition}

\begin{theorem} \label{main1} Let
$$
\mathfrak S=\prod_{n=1}^\infty S_n^{f(n)}
$$
be a Cartesian product of an infinite  family $\{S_n\}$ of finite simple groups
such that $l(S_n) \rightarrow \infty$. If $\mathfrak S$ is topologically finitely 
generated then it is a profinite completion.
\end{theorem}

We say that the family $\{S_n\}$ of simple groups has 
\emph{essentially unbounded rank} if $l(S_n) \rightarrow \infty$.
Note that if $S$ is a classical simple group of rank $r(S)$ then the ratio $l(S)/r(S)$ is bounded. 
Thus $l(S_n) \rightarrow \infty$ is equivalent to
the ranks $r(S_n)$ tending to infinity.
Example~\ref{ex} shows that some restriction of this form in Theorem~\ref{main1} is necessary.

It is easy to see that the product $\mathfrak{S}$ in Theorem~\ref{main1} is 
finitely generated iff there is $c$
such that $f(n)<|S_n|^c$ for all
$n \in \mathbb N$. We shall prove that then $\mathfrak{G}$ has a 
\emph{frame} subgroup $\Gamma$
as defined below. Moreover, $\Gamma$ can be taken to be 
$22(c+1)$-generated. 

\begin{remark} Actually our proof shows the following more general result:

Let $d \in \mathbb{N}$.
Suppose that $\{T_n\}$ is a family of finite groups.
Assume that each $T_n$ is generated by the images of $d$ homomorphisms 
from the alternating group $\Alt(n)$ to $T_n$.
Then the product $\prod_{n=1}^\infty T_n$ is a profinite completion of a 
finitely generated group.
\end{remark}

Suppose that a Cartesian product $\mathfrak{G}$ of alternating groups 
is topologically generated by $d$ elements.
In view of Theorem~\ref{main} the following natural question presents itself:

\begin{question}
Is $\mathfrak{G}$ the completion of a $d$-generated residually finite group?
\end{question}
This seems unlikely but we don't have a counterexample.
\newpage

Our interest in Cartesian products as profinite completions was motivated by our search for
possible connections between property $\tau$ and subgroup growth.
In this setting we ask the following natural (and seemingly difficult) question:
\begin{question}
Is $\G$ the profinite completion of a finitely generated group which has property $\tau$?
\end{question}
At least some Cartesian products $\G$ are, namely those 
$\G= \prod_{n=1}^\infty \Alt(n)^{f(n)}$ with $f(n) < n^{\log n}$. See
\cite{martin3} to which this note is a companion. 
However at this time we don't even know whether
in general $\G$ always has a dense finitely generated subgroup with $\tau$.

\bigskip

Our proof is based on the idea of \emph{frame} subgroups of Cartesian products. 
These are defined in Section~\ref{frames}.
Theorem~\ref{main1} is then reduced in Section \ref{easy} to proving the existence 
of a single frame subgroup (namely Theorem~\ref{beginning}).
This is the main technical difficulty and it is done in Section~\ref{the hard work} 
using some results of the first author in \cite{martin}.

\subsection*{Notation}
\begin{itemize}

\item $\Alt(n)$ (resp. $\Sym(n)$) is the alternating
(resp. symmetric) group on $n$ letters,

\item $[a,b]$ is the set of integers between $a<b$,

\item The elements of a Cartesian product 
$\prod_{n \in I} S_n$ of groups are denoted by $(a_n)_{n \in I}$
or just by $(a_n)_n$ when there is no possibility of confusion (each $a_n \in S_n$).
\end{itemize}

\section{Frame subgroups}
\label{frames}

Let $\mathfrak{S} = \prod_{n=1}^\infty S_n$ 
be a Cartesian product of finite groups.

 \begin{definition} \label{deframe}
A finitely generated subgroup $G < \mathfrak{S}$ is
\emph{a frame} for $\mathfrak{S}$ if the following hold: 

\medskip

(a) $G$ contains $\bigoplus_{n=1}^\infty S_n$. 

\medskip

(b) The natural surjection $\widehat G \rightarrow \mathfrak{S}$ is an isomorphism.
\end{definition}

One can think of condition (a) as saying that $G$ is a good approximation of 
$\mathfrak{S}$ from 'within'
while condition (b) says that $G$ approximates very well $\mathfrak{S}$ from 'above'. 

\medskip

More precisely, for a finite subset $V \subset \mathbb N$ of integers define the 
\emph{$V$-principal congruence subgroup}
 $G_V$ to be the kernel of the projection of $G$ onto $\prod_{n \in V} S_n$.
Let $G(V)$ be the projection of $G$ onto $\mathfrak{S} (V):= \prod_{n \not \in V} S_n$.
The $m$-th principal congruence subgroup $G_m$ is just 
$G_{ \{1,\ldots ,m\}}$ and $G (m):=  G (\{1,\ldots ,m\})$.
\medskip

Part (a) of the above definition is now equivalent to
$$
G = \left(\prod_{n\in V} S_n \right) \times  G_V,
\quad
G_V = G \cap \mathfrak{S} (V).
$$
Therefore the congruence subgroup $G_V$ can be identified with the projection $G(V)$. \medskip

On the other hand, part (b) of Definition~\ref{deframe} says that 
the profinite topology of $G$ is the same as its congruence topology:
Every subgroup of finite index in $G$ contains a congruence subgroup 
$G_m$ for some $m \in \mathbb{N}$.

\medskip

We stress that the existence of even a single example of 
a frame is far from obvious at this stage.
\medskip

The following Lemma allows us to find many frame groups provided we already know at least one:

\begin{lemma}
\label{glue}
Let $A_n,B_n <C_n$, ($n \in \mathbb{N}$) be finite groups with
$C_n= \langle A_n,B_n \rangle$. Suppose that $X$ (resp. $Y$) is a frame subgroup
of  the product $\mathfrak A= \prod_{n=1}^\infty A_n$ (resp. $ \mathfrak B =\prod_{n=1}^\infty B_n$).  
Each of $X$ and $Y$ can
be considered as a subgroup of $\mathfrak{C}:= \prod_{n=1}^\infty C_n$ in the natural way. 
Then the group
$$
Z=\langle X,Y \rangle < \mathfrak{C}
$$
is a frame in $\mathfrak{C}$.
\end{lemma}

\textbf{Proof:} It is clear that $Z$ contains the direct product of $C_n$.
Suppose that $N$ is a subgroup of finite index in $Z$. 
By hypothesis $N$ contains the $m$-th principal congruence subgroups
$X_m$, $Y_m$ for some $m$. Identifying them with $X(m),Y(m)$ in $\mathfrak{C} (m)$ we see that
$N$ contains $\langle X(m),Y(m) \rangle = Z(m)$, which under our identification is $Z_m$.
$\square$

It is clear that for every $c>2$ Lemma~\ref{glue} can be generalized 
by induction from a pair $X,Y$ to
$c$ frame subgroups $X_1,\ldots ,X_c$ satisfying an analogous condition.
\medskip

For possible future use we prove the following extension:
\begin{lemma}
\label{extension}

Let $A_n,M_n,B_n$, ($n= 1,2, \ldots $) be finite groups such that
$M_n= A_n\ltimes B_n$.
For each $n$ let $b_{n,1},\dots,b_{n,k}$ be elements in  $B_n$, such that $M_n$ is generated
by $A_n$ and $[A_n,b_{n,s}]$ ($1 \leq s \leq k$).

Suppose that $X=\langle x_1,\dots,x_m\rangle$ is a frame subgroup
of  the product $\prod_{n=1}^\infty A_n$.  Then $X$ can
be considered as a subgroup of $\mathfrak{M}:= \prod_{n=1}^\infty M_n$ in the natural way.
Define $b_s = (b_{n,s})_n$ for $s=1,\ldots k$.
Then the group
$$
Z=\left\langle X,  [x_j,b_s] \ | \ 1\leq j \leq m,\ 1\leq s \leq k \right\rangle 
< \mathfrak{M}
$$
is a frame in $\mathfrak{M}$.
\end{lemma}

\textbf{Proof:} The condition that $M_n$ is generated by $A_n$ and $[A_n,b_{n,s}]$, 
together with $X \geq \oplus A_n$, easily gives that
$Z \geq \oplus M_n$.

Suppose now that $H$ is a normal subgroup of finite index in $Z$. 
By hypothesis $H$ contains the $l$-th principal congruence subgroup
$X_l$ for some $l$. Identifying it with $X(l)$ in $\mathfrak{M} (l)$ we see that
$N$ contains $\langle X(l), [X(l),n_s] \mid \ s=1,\ldots ,k \rangle=Z(l)$, 
which under our identification is $Z_l$.
$\square$

\section{Reductions} \label{easy}

The following is a corollary of \cite{TW} and \cite{lsh}, Theorem 1.2:
\begin{theorem} \label{useful}
For every $m \geq 5$ there is an integer $r=r(m)$ such that if $S$ is a 
finite simple group with $l(S)>r$ then $S$
is generated by two 
subgroups isomorphic to $\Alt(m)$.
\end{theorem}

Incidentally this raises the following basic
\begin{question}
Is it true that if a finite simple group $S$ contains a copy of 
$\Alt(m)$, $m \geq 5$ then $S$ is generated by two copies of
$\Alt(m)$?
\end{question}

\bigskip

Our starting point for Theorem \ref{main1} is the following
\begin{theorem}
\label{beginning}
For every odd prime $p$, there exists a 10-generated group
$G_1$ which is a frame for the Cartesian product
$$
\prod_{n=3}^\infty \Alt( u_{n,p}),
$$
where $u_{n,p}=(p^{3n}-1)(p-1)^{-1}$.
\end{theorem}

This is proved in Section~\ref{the hard work}.
Assuming that we complete the proof of Theorem~\ref{main1} in two steps as follows:
\subsubsection*{Step 1}
\begin{proposition}
\label{st1}
Given any sequence $\{ S_n\}$ of distinct
finite simple groups of essentially infinite rank, there is a 22-generated group 
$G_2$ which is frame for
the Cartesian product
$$
\G= \prod_{n=1}^\infty S_n .
$$
\end{proposition}

\textbf{Proof:} Without loss of generality we may assume that $S_n$ 
are numbered so that $l(S_1) \leq l(S_2) \leq \ldots $.

Assume first that $l(S_1) \geq r(3^9)$, where $r(m)$ is the number from 
Theorem~\ref{useful}. \medskip

Define a function $h:\mathbb{N} \rightarrow \mathbb{N}$ inductively by

$h(1)=1$,

$h(k)$ is the smallest $n>h(k-1)$ such that $l(S_n) \geq r( u_{k})$, where 
$u_k:=(3^{3(k+2)}-1)/2=u_{k+2,3}$.
\medskip

The existence of such $h(k)$ follows
from the fact that $l(S_i) \rightarrow \infty$.

Set
$$
B_k= \prod_{ h(k) \leq n  < h(k+1)} S_n.
$$
Then $\G= \prod_{k=1}^\infty B_k$.

For every $n \in [h(k),h(k+1))$ we have that $l(S_n) \geq r(u_k)$ and 
therefore by Theorem~\ref{useful} the group $S_n$ is generated by two copies of $\Alt(u_{k})$.
In other words we have two embeddings
$$
f_{n,j}: \ \Alt(u_k) \rightarrow S_n, \quad i=1,2,
$$
such that
$S_n=\langle f_{n,1} (\Alt(u_k)) ,f_{n,2} (\Alt(u_k)) \rangle$.

Now $B_k$ is generated by two copies of $\Alt(u_k)$ as follows:

For $j=1,2$ define $P_{k,j}$ to be the image of
$$
\Alt(u_k) \ni a \mapsto (f_{n,j}(a))_n  \in B_k, \quad (h(k) \leq n <h(k+1)).
$$

Then $\langle P_{k,1}, P_{k,2} \rangle$ is the whole of $B_k$ 
because it is a subdirect product of distinct simple groups.

\medskip

Now by Theorem~\ref{beginning} there are two embeddings, $t_j$ ($j=1,2$) of
$G_1$ into $\prod_{k=1}^\infty P_{k,j}$ such that the images $t_j(G_1)$ are frame subgroups.
Lemma~\ref{glue} gives that these two copies of $G_1$ embedded in $\prod_{k=1}^\infty B_k=\G$
via the $t_j$ generate a frame subgroup $G_2$ for $\G$.

In general we can
write $\G$ as $\G= K \times \G'$ where $K$ is a finite product of distinct simple groups and
$\G'$ is a Cartesian product of simple groups $S$ with $l(S)>r(3^9)$. 
Taking a 20-generated frame $G'$ in $\G'$ together with 2 generators $a,b \in K$
gives the frame $G=\langle a,b, G' \rangle$ in $\G$.
$\square$

\subsubsection*{Step 2}
Recall the following well known
\begin{proposition}
\label{multiple}
For any finite simple group $S$ and any integers 
$c,m \in \mathbb N$  such that $m\leq  |\Aut(S)|^c$
there exist $c+1$ embeddings
$$
f_i: S \rightarrow  S^{m} \quad (i=1,\ldots ,c+1)
$$
such that $S^m=\langle f_1(S), \ldots ,f_{c+1}(S) \rangle$.
\end{proposition}

\textbf{Proof:}
It is sufficient to prove the Proposition for the case when $m=|\Aut (S)|^c$.

Set $N=|\Aut(S)|$.
Identify $S^{N^c}$ with the sequences $(a_\mathbf{t})_\mathbf{t}$ labelled by the
$c$-tuples $\mathbf{t}=(t_1,\ldots , t_c) \in \Aut (S) ^c$ of elements of $\Aut (S)$.

Now define $c+1$ embeddings $f_i : S \rightarrow S^{N^c}$ by
$f_0: \ a\in S \mapsto (a)_\mathbf{t}$ for $i=0$ and for $i=1,2,\ldots c$ define
$$
f_i: \quad a \mapsto (a^{t_i})_\mathbf{t}, \quad  a\in S.
$$
We claim that the $(c+1)$ subgroups $f_i(S)$ generate $D$:

By \cite{cameron} Exercise 4.3, if $M < \prod_{k=1}^n{S_k}$ is a subdirect product
of isomorphic non-abelian simple groups $S_k$, 
then there exist indices $1\leq i<j \leq n$ such that the projection of $M$ onto
$S_i \times S_j$ is the diagonal subgroup $\{(a,f(a))\ | \ a \in S_i \}$, 
where $f: S_i \rightarrow S_j$ is an isomorphism.

We have chosen the emebeddings $f_i$ such that this is impossible.
$\square$
\medskip

We also need

\begin{proposition}
\label{generators}
For every $c \in \mathbb{N}$, every finite simple group $S$ and $m>|S|^c$, the direct product
$D=  S^m$ is not generated by $c$ elements.
\end{proposition}
\textbf{Proof:} Suppose that $g_i=(g_{i,k})_{k=1}^m \in D$, $(i=1,\ldots c)$ 
are any $c$ elements of $D$. Since $m>|S|^c$
there are two coordinates $k \not =k' \in \{1,\ldots m\}$ such that 
$g_{i,k}=g_{i,k'}$ for all $i\in [1,c]$. Therefore $g_1,\ldots,g_c$ cannot
generate $D$.

In fact a similar proof shows that $D$ is not $c$-generated for $m>\frac{|S|^c}{|\Aut(S)|}$.
$\square$

\bigskip

Now we can easily finish the proof of Theorem~\ref{main1}:
\medskip

Suppose $\G=\prod_n S_n^{f(n)}$ is topologically finitely generated.

By Proposition~\ref{st1} there is a $22$-generated group $G_2$ which is a frame in $\prod_n S_n$.
Since $\G$ is assumed to be finitely generated, by
Proposition \ref{generators} there is $c\in \mathbb N$ such that
$$
f(n)<|S_n|^c \leq |\Aut(S_n)|^c,\,\,\,\,\,\mbox{for all }n.
$$
Now by
Proposition~\ref{multiple}, for each $n=1,2,\ldots $ we find 
$(c+1)$ copies of $S_n$ inside the product
$S_n^{f(n)}$, which together generate it. 
By an application of Lemma~\ref{glue} we deduce that there
exists a frame subgroup $\Gamma < \G$ generated by $(c+1)$
copies of $G_2$, and Theorem~\ref{main1} is proved.

$\square$ \medskip

We conclude this Section with
\begin{question}
Does there exists a function $f(d)$ with the following property:

Suppose that $n \in \mathbb N$ and $G$ is a $d$-generated finite group, which 
can also be generated by several subgroups isomorphic to $\Alt(n)$. Then 
there exist $f(d)$ homomorphisms 
$i_s: \Alt(n) \to G$, such that
$$
G = \langle \mathrm{im\, } i_s \mid s=1,\dots, f(d) \rangle.
$$
\end{question}

If true this will imply the following generalization of Theorem \ref{main1}: 

A finitely generated Cartesian product $\prod G_n$ is a profinite completion, 
provided that each $G_n$ is generated by subgroups isomorphic to $\Alt(n)$.

\section{Theorem~\ref{beginning}}
\label{the hard work}

\subsection{Rings and $\EL_3$}
In this section we present a construction from~\cite{martin} which will be 
necessary for Theorem~\ref{beginning}.
By $M_n(A)$ we denote the ring of $n\times n$ matrices over a ring $A$. The
elementary matrix in $M_n(\F_p)$ with entry $r\in A$ in position $(i,j)$ and $0$ 
elsewhere is $e_{i,j}(a)$. We write
$e_{i,j}:=e_{i,j}(1)$.

Let $R$ be the subring of $\prod_{n=3}^\infty M_n(\F_p)$ generated by the
five elements
$1= (\textrm{Id}_n)$,
$\mathbf{a}=(a_n),\ \mathbf{a}^{-1}, \mathbf{b}=
(b_n)$ and $\mathbf{c}=(c_n)$, defined as follows:
$$
a_n=e_{1,2}+e_{2,3}+\ldots +e_{n,1}, \quad b_n=e_{1,2}, \quad c_n=e_{2,1}, \quad i \geq 2 .
$$
Note that $R$ contains the elements
$$
\mathbf d := [\mathbf b,\mathbf a^{-1} \mathbf b \mathbf{a}] = (e_{1,3})_n ,
\quad \mathbf e: = [\mathbf a^{-1} \mathbf c \mathbf{a},\mathbf c] = (e_{3,1})_n,
$$
which together with $\mathbf{b}$ and $\mathbf{c}$ and $1$ generate a subring of 
$R$ isomorphic to $M_3(\F_p)$ sitting diagonally
in the top left corner of each factor.

\begin{proposition}[\cite{martin}]
\label{pro-R}
The ring $R$ contains the direct sum
$$
\bigoplus_{n=3}^\infty M_n(\F_p).
$$
The profinite completion $\hat R$ of $R$ is
$$
\widehat R= \widehat{\F_p[t,t^{-1}]} \bigoplus  \prod_{n=3}^\infty M_n(\F_p),
$$
\end{proposition}
\textbf{Proof:} (sketch)

Let us observe that the commutator $[b_n,a_n^{-k}b_na_n^k]$ is non-zero iff 
$n$ divides $k-1$ or $k+1$.
This, together with the simplicity of the ring $M_n(\F_p)$ implies that the ring generated by
$\mathbf{a}$, $\mathbf{a}^{-1}$ and $\mathbf{b}$ contains the infinite direct sum.

Let $I$ be an ideal of finite index in $R$, which contains $\mathbf{a}^N-1$.
Let $\tilde R = R/J$, where $J= I + \oplus_{n=3}^{2N} M_n(\F_p)$.
Then in $\tilde R$ we have
$$
\tilde{\mathbf{d}} =
   \left[\tilde{\mathbf{b}},\tilde{\mathbf{a}}^{-1} \tilde{\mathbf{b}} \tilde{\mathbf{a}}\right] =
\left[\tilde{\mathbf{b}},\tilde{\mathbf{a}}^{-N-1} \tilde{\mathbf{b}} \tilde{\mathbf{a}}^{N+1}\right] 
= 0
$$
because
$$
\left[{\mathbf{b}},{\mathbf{a}}^{-N-1} {\mathbf{b}} {\mathbf{a}}^{N+1}\right] 
\subset \bigoplus_{n=3}^{2N} M_n(\F_p)\subset J.
$$
Similar proof gives also that $\tilde{\mathbf{e}} =0$. The relations
$$
\mathbf{b} = \left[\mathbf{d},\mathbf{a}^{-1} \mathbf{c}\mathbf{a}\right],
\quad
\mathbf{c} = \left[\mathbf{a}^{-1} \mathbf{b} \mathbf{a},\mathbf{e}\right]
$$
give that $\tilde R$ is an image of $\F_p C_N=\F_p[t,t^{-1}]/\langle t^N=1 \rangle$ 
because it is generated by $\mathbf{a}$ and $\mathbf{a}^{-1}$,
which implies that the profinite completion of $R$ has the desired form.
$\square$

The profinite ring $U:=\widehat{\F_p[t,t^{-1}]}$ 
is the inverse limit of the group algebras $\F_p C_n$
of the cyclic groups $C_n$ for $n \in \mathbb N$. If $\mathrm{proj}_U$ is the projection of
$\hat R$ onto $U$ then the above argument gives that:

$$ \mathrm{proj}_U (\mathbf{a}^\pm)= t^\pm ,\quad \mathrm{proj}_U
(\mathbf{b})= \mathrm{proj}_U(\mathbf{c})=0.
$$

\begin{definition}
Let $G_0:=\EL_3(R)$ be the subgroup of $M_3(R)$ generated by the elementary matrices $e_{i,j}(r)$ for
$1\leq i\not =j \leq 3$ and $r\in R$.
\end{definition}
Then $G_0$ is generated by the following set of ten elements:
$$
\{e_{i,j} \ \mid \ 1\leq i\not = j \leq 3\} \ \bigcup \ I 
: = \left\{ e_{1,2}(x) \ \mid  \ x \in \{\mathbf{a,a^{-1},b,c} \}\right\}.
$$

In fact it is easy to see that $G_0$ is generated by the $5$ elements

$$
\left\{ e_{1,2}(\mathbf{a}),e_{2,3}(\mathbf{a^{-1}}),e_{1,2}(\mathbf{b}),e_{1,2}(\mathbf{c}),g\right\},
$$
where $g$ is a matrix in $\SL_3(\F_p)$ such that $\SL_3(\F_p)=\langle e_{1,3}, g \rangle$.
\medskip

Since $R$ is a subring of $\prod_{n=3}^\infty M_n(\F_p)$ we can consider $G_0$ as a subgroup of
$$
\prod_{n=3}^\infty \EL_3(M_n(\F_p))= \prod_{n=3}^\infty \SL_{3n}(\F_p).
$$

\begin{proposition}[\cite{martin}]
\label{g0}
The group $G_0$ contains
$$
\bigoplus_{n=3}^\infty \SL_{3n}(\F_p).
$$

The profinite completion of $G_0$ is
$$
\EL_3(U) \oplus \prod_{n=3}^\infty \EL_3(M_n(\F_p))=
\lim_{\substack{ \longleftarrow \\ n\in \mathbb{N} }} \SL_3 (\F_p C_n)
\oplus \prod_{n=3} ^\infty \SL_{3n}(\F_p).
$$

We can say more. Suppose $\bar G_0=G_0/H$  is a finite image of $G_0$. Then:
\medskip

(i) There exists $N < [G_0:H]$ such that we have the following diagram:
$$
\begin{array}{ccccc}
               &          & G_0 &          &   \\
               & \swarrow &     & \searrow & \\
\SL_3(\F_pC_N) \!\!\!\!\!\!\!\!\!\!\!\!&          &     &          & \!\!\!\!\!\!\!\!\!\!\!\! G_0/H \\
               & \searrow &     & \swarrow & \\
               &          & \!\!\!\!\!\!\!\!\!\!\!\!G_0/\widetilde H =\widetilde{G}_0 \!\!\!\!\!\!\!\!\!\!\!\!
                                &       &
\end{array}
$$
where $\widetilde H =H \cdot \oplus_{n=3}^N \SL_{3n}(\F_p)$. The map $G_0 \to \SL_3(\F_pC_n)$ comes
from the projection $R \to \F_pC_N$, sending $\mathbf{b}$ and $\mathbf{c}$ to $0$. 
\medskip

(ii) If $\pi$ is the map $G_0 \to \widetilde{G}_0$ then for each pair of indices 
$1\leq i \not = j \leq 3$ we have
$
\pi (e_{i,j}(\mathbf{b}))= \pi (e_{i,j}(\mathbf{c})) =1$,
and moreover, $\widetilde{G}_0$ is the normal closure of any of its elements $\pi (e_{i,j})$. 
\medskip

(iii) $G_0/H = \widetilde{G}_0 \times M$, where
$M$ is a central quotient of $\oplus_{n=3}^N \SL_{3n}(\F_p) <G_0$.

\end{proposition}
\textbf{Proof:}
Using Proposition~\ref{pro-R} we can see that
the group $G_0=\EL_3(R)$ has the congruence subgroup property (CSP) because
$\EL_3(\F_p[t,t^{-1}])$ has CSP and $K_2(M_n(\F_p))$ is trivial for any $n$.
Therefore the profinite completion of $G_0$ is the same as the congruence completion which is
$$
\EL_3(U) \oplus \prod_{n=3}^\infty \EL_3(M_n(\F_p)).
$$
This is easily gives (i) and (ii). Part (iii) follows from the fact that 
$G_0$ contains $\bigoplus_{n=3}^\infty \SL_{3n}(\F_p)$.
$\square$

\subsection{The construction}

The group $G_0$ is our first `approximation' to a frame subgroup. 
It remains to modify it and suppress the
`bad' factor $\EL_3(U)$ in $\hat G_0$.

Let
$$
D:=\left \langle e_{1,3}(\mathbf{b}), e_{3,1}(\mathbf{c}), e_{1,3}(\mathbf{d}), e_{3,1}(\mathbf{e})
\right \rangle <G_0.
$$
Then $D$ is isomorphic to a copy of
$$
\SL_3(\F_p) \simeq D <G_0< \prod_{n=3}^\infty \SL_{3n}(\F_p)
$$
diagonally embedded in positions
$(a,b) \in \{1,2n+2,2n+3\}^2$ in each factor $\SL_{3n}(\F_p)$.  
\medskip

We remark that from the Chevalley commutator relations for elementary matrices it follows that
$e_{1,2}(\mathbf{d})$ and $e_{1,2}(\mathbf{e})$ are each in the normal closure of
$\langle e_{1,2}(\mathbf{b}), e_{2,1}(\mathbf{c}) \rangle$ in $G_0$.
Therefore, the projection of $D$ in the `bad' factor 
$\EL_3(U)$ of $\widehat{G}_0$ is trivial, by Proposition~\ref{g0} (ii).
\medskip

Define $q:=(q_n) \in D$ where $q_n \in \SL_{3n} (\F_p)$
is the diagonal element with entries all 1 except
$-1$ in positions $2n+2$ and $2n+3$ on the diagonal of $\SL_{3n} (\F_p)$.
It is clear that $q$ commutes with the element
$v:=e_{1,2}=(e_{1,2} (\mathrm{Id}_n)) \in G_0$.
\medskip

Recall that $u_{n,p}=(p^{3n}-1)/(p-1)$. We will use the
following crucial Lemma:
\begin{lemma}
\label{crux}
For every $n \geq 2$ there exist two homomorphisms $i'=i'_n, \ i''=i''_n$
of $K:=\SL_{3n}(\F_p)$ into $\Alt(u_{n,p})$ such that
\begin{equation}\label{rel}
\begin{array}{c}
\displaystyle
i''(q)^{-1}i'(v)i''(q)=i'(v)^{-1},
\quad
i'(q)^{-1}i''(v)i'(q)=i''(v)^{-1}
\\
\\
\displaystyle
\mathrm{and }\quad
\langle i'(K),i''(K) \rangle = \Alt(u_{n,p}).
\end{array}
\end{equation}
\end{lemma}

Assuming this we can prove Theorem~\ref{beginning}: 
\medskip

Set $S_n:=\Alt(u_{n,p})$.
Consider the two homomorphisms
$$
E',E'': \ \prod_{n=3}^\infty \SL_{3n}(\F_p) \rightarrow \mathfrak A:=\prod_{n=3}^\infty S_n,
$$
given by $(g_n)_n \mapsto (i'_n(g_n))_n$ and $(g_n)_n \mapsto (i''_n(g_n))_n$.  
\medskip

Set $G'_0=E'(G_0), G''_0=E''(G_0)$. For an element $a \in G_0$ we shall write
$a'$ (resp. $a''$)  for $E'(a) \in G'_0$ (resp. $E''(a) \in G''_0$) 
and do the same for subgroups of $G_0$. \medskip

\begin{definition}
Set $G=G_1:=\langle G'_0, G''_0 \rangle < \mathfrak A$.
\end{definition}

We shall prove that $G_1$ is a frame subgroup for $\mathfrak A$. 

The group $G_1$ contains $\bigoplus_{n=3}^\infty S_n$, because each copy of $G_0$ contains the
infinite direct sum of the images of $\SL_{3n}(\F_p)$ in $S_n$.
\medskip


Suppose now that $G_1/H$ is a finite quotient of $G=G_1$.
By Proposition~\ref{g0} (i) applied to $G'_0$ and $G''_0$, 
there exists an integer $N < [G_1:H]$ such that if
$$
\widetilde{H}  = H . \bigoplus_{n=3}^N S_n \quad \widetilde{H'} = \widetilde H \cap G'_0 , 
\quad
\widetilde{H''}= \widetilde H \cap G''_0.
$$
then
$\widetilde{G} = G_1/\widetilde{H}$ is generated by
$\widetilde{G'_0} = G'_0/\widetilde{H'}$ and 
$\widetilde{G''_0} = G''_0/\widetilde{H''}$ which are images of $\SL_3(\F_pC_N)$. 

By Proposition~\ref{g0} (ii) the images of 
$e_{i,j}(\mathbf{b})$ and $e_{i,j}(\mathbf{c})$ in 
$\widetilde{G'_0}$ and $\widetilde{G''_0}$ are trivial.
In particular the group $D$ is in the kernel of both 
$\pi'= \pi \circ E'$ and $\pi''= \pi \circ E''$,
where $\pi$ is the projection $G \rightarrow \widetilde{G}$.
This implies that $\pi'(q) = \pi''(q) =1$.

By (\ref{rel}) the following relations hold in $G$:
$$
E'(q)^{-1} E''(v) E'(q) = E''(v)^{-1}
\quad\quad
E''(q)^{-1} E'(v) E''(q) = E'(v)^{-1}.
$$
Therefore we have
$$
\pi'(q)^{-1} \pi''(v) \pi'(q) = \pi''(v)^{-1}
\quad\quad
\pi''(q)^{-1} \pi'(v) \pi''(q) = \pi'(v)^{-1}.
$$
These, together with $\pi'(q)=\pi''(q)=1$ and $v^p=1$ imply that $\pi'(v)=\pi''(v)=1$.
The last two elements normally generate $\widetilde{G'}$ and $\widetilde{G''}$, 
therefore these groups are trivial. This shows that $\widetilde{G}=
\{1\}$ , i.e., $G = H .\oplus_{n=3}^N S_n$. 
Now $G> \oplus_{n=3}^N S_n$, hence $H$ contains the congruence subgroup $G_N$ 
and so $G$ is a frame.
$\square$

Theorem~\ref{beginning} is now proved modulo Lemma~\ref{crux}.

\subsection{Proof of Lemma~\ref{crux}}

Consider the action of $\SL_{3n}(\F_p)$ on the set $P$ of 
$u_{p,n}$ points of the projective plane $\mathbb{P}(V)$, where $V=\F_p^{3n}$.

This gives the homomorphism $i': \ \SL_{3n}(\F_p) \rightarrow \Alt(u_{p,n})=:S$.
The element $q$ is a diagonal element with two eigenvalues
$-1$ and so $q'=i'(q)\in S$ is an involution which fixes $(p^{3n-2}+p^2-2)/(p-1)$ points of $P$.
The element $v'=i'(v)=i'(e_{1,2})$ has order $p$ and we see that it moves a set 
$M \subseteq P$ consisting of $mp$ points in
$m=(p^{3n}-p^{2n})p^{-1}(p-1)^{-1}$ $p$-cycles. Of these $p$-cycles
$m/p^2$ are pointwise fixed by $q'$ and the rest are swapped by $q'$. 
Recall that $v'$ and $q'$ commute.

\bigskip



The essential point now is that exactly $1/p^2$ of the $p$-cycles of $v'$ are
fixed by the involution $q'$ and the rest are 'swapped'.
Without loss of generality we may assume that the set $M$ of points moved by 
$v'$ is arranged in $N=m/p^2$ cubes of size $p$ each, say
$$
M= \bigcup_{i=1}^{N} B_i, \quad B_i = \{b_i(j,k,m) \ | \ j,k,m \in [0,p) \}.
$$

We may assume that $v'$ acts on $M$ by a shift of the first coordinate
$b_i(j,k,m) \mapsto b_i(j+1,k,m)$, while $q_i$ restricted to $M$
acts by
$$
b_i(j,k,m) \mapsto b_i(j,-k,-m)   \quad \forall i \in [1,N],\ j,k,m \in [0,p).
$$
The indices $j,k,m$ are taken mod $p$ everywhere.

\bigskip

Let $\sigma \in \Sym (P)$ be defined on $M$ as the cyclic shift
$$
b_i(j,k,m) \mapsto b_i(k,m,j)
$$
and on $P \setminus M$ as the the identity.

Consider $q''=(q')^\sigma$ and $v''=v'^\sigma$. 
Then $q''$ acts by a shift of the second coordinate,
while $v''$ acts as the reflection in the 1-st and 3-rd coordinates. 
An easy calculation now shows that
the conditions (\ref{rel}) hold.

So we can define
$$
 i'':\  \SL_{3n}(\F_p) \rightarrow S, \quad i''(x)= (i'(x))^\sigma
$$
provided $\langle (K')^\sigma, K'\rangle =S$. \medskip

It is enough to show that $(K')^\sigma \not = K'$:

By~\cite{kantor} if $K'$ is contained in a proper subgroup $H < S=\Alt(u_{p,n})$ then $H \leq N_S(K')$. 
Assuming $K''=i''(K)=(K')^\sigma \not = K'$ then $K''$ certainly
does not normalize $K'$ and so $\langle K',(K')^\sigma \rangle =S$. \medskip

Now suppose that $\sigma$ defined above normalizes $K'$. 
Replace $\sigma$ with $\sigma' =a \sigma $ where $a$ is a transposition which centralizes
both $q'$ and $v'$. For example we can take $a$ to be any transposition moving points 
of $P$ fixed by both $q'$ and $v'$: there are
$(p^{2n-2}+p^2-2)/(p-1)$ of them. Then 
$(v')^{\sigma'}=(v')^\sigma$, $(q')^{\sigma'}=(q')^\sigma$ and if $\sigma'$ 
normalizes $K'$ then so does
$a$. However, in that case $\langle K',a \rangle$ is a proper primitive 
subgroup of $\Sym (P)$ containing a transposition, which is impossible.
$\square$

\begin{acknowledgement} 
Part of this research was carried out during a visit of the first author to Oxford in January 2006. 
We gratefully acknowledge the support of the London Mathematical Society with a Scheme 4 grant.
\end{acknowledgement}

\texttt{\\ Martin Kassabov, \\ Department of Mathematics,\\ Cornell University, Ithaca, NY 14853, USA. \\ kassabov@math.cornell.edu 
\medskip
 \\ Nikolay Nikolov, \\ New College, Oxford, OX1 3BN, UK. \\ nikolay.nikolov@new.ox.ac.uk }
\end{document}